\documentclass{gtart}
\usepackage{amsthm,amssymb}

\input gtmonout
\volumenumber{2}
\volumeyear{1999}
\volumename{Proceedings of the Kirbyfest}
\pagenumbers{113}{123}
\papernumber{6}
\received{12 March 1999}\revised{19 July 1999}
\published{17 November 1999}

\newtheorem{thm}{Theorem}[section]

\newtheorem{prop}[thm]{Proposition}

\input rlepsf

\begin{document}

\title{$\mathbf{Z_2}$--Systolic-Freedom}
\covertitle{$\noexpand\bf Z_2$--Systolic-Freedom}
\asciititle{Z_2-Systolic-Freedom}
\author{Michael H Freedman}

\address{Microsoft Corporation, One Microsoft Way\\
Redmond, WA 98052-6399, USA}
\email{michaelf@microsoft.com}

\begin{abstract}  We give the first example of systolic freedom
over torsion coefficients. The phenomenon is a bit unexpected
(contrary to a conjecture of Gromov's) and more delicate than
systolic freedom over the integers.

{\em Dedicated to Rob Kirby, a lover of Mathematics and
other wild places. Thank you for your inspiration.}

\end{abstract}

\asciiabstract{We give the first example of systolic freedom
over torsion coefficients. The phenomenon is a bit unexpected
(contrary to a conjecture of Gromov's) and more delicate than
systolic freedom over the integers.}

\maketitle

\setcounter{section}{-1}
\section{Introduction} For closed Riemannian surfaces, whose
topology is different from the $2$--sphere, $$ A \geq \frac{2}
{\pi} L^2\eqno(0.1) $$ where
$A$ is area and $L$ is the length of the shortest essential
loop. The boundary case is a round projective plane.  See
\cite{KS} and
\cite{G1} for a discussion.  For closed manifolds of higher
dimensions, such $\lq\lq$systolic inequalities" have been the
focus of much research and many interesting counter-examples
exist \cite{BK}, \cite{G2}, and \cite{G3}.

We recall some definitions:\\ Let $M$ be a closed Riemannian
manifold of dimension $n$ and let $0 \leq p,\, q \leq
n,\linebreak p+q=n$. $$ systole_k(M)=\inf\,
area_k[\alpha]\eqno(0.2) $$ where the infimum is taken over all
smooth oriented $k$--cycles $\alpha$ with $[\alpha]\neq 0 \in
H_k (M;Z).$ $$ Z_2-systole_k(M)=
\inf\,area_k(\alpha) \eqno(0.3) $$ where the infimum is taken
over unoriented $k$--cycles $\alpha$, $[\alpha]\neq 0 \in
H_k(M;Z_2).$
$$ stable-systole_k(M)=\inf\,\,stable- area_k[\alpha]\eqno(0.4)
$$ where $[\alpha]\neq 0 \in H_k(M;Z)/$torsion and $$
stable-area_k
\alpha=\inf\frac {1}{i}\Big(\inf_{i[\alpha]}area_k (\beta)\Big)
$$ where $i=1,2,3,\dots$ in the inner infimum is over oriented
cycles
$\beta$ representing $i[\alpha].$

Gromov proved (see \cite{H} for discussion and generalizations)
that
$\lq\lq$stable systolic rigidity" holds for any product of
spheres
$S^p \times S^q =: M^n$, that is there a constant $c(n)$ so that
for any Riemannian metric on $M^n = S^p \times S^q ,p+q=n$: $$
vol(M)=stable-systole_n(M)\geq c \cdot
stable-systole_p(M)\cdot stable-systole_q(M)\eqno(0.5) $$
Surprisingly, he also discovered that the corresponding
unstable statement is false:\\ Let
$M_r=S_r^3 \times \mathbb{R}/(\theta, t)\equiv(\sqrt{r} \circ
\theta, t+1)$, where $S_r^3$ is the $3$--sphere of radius $r$ and
the identification matches a point with its $\sqrt{r}$--rotation
along Hopf fibers displaced one unit in the real coordinate. For
this $r$--family of metrics on $S^3 \times S^1$, we have
$\lq\lq(3,1)$--systolic freedom" $$
\frac{systole_4(M_r)}{systole_3(M_r)\,systole_1(M_r)}=
\frac{\mathcal{O}(r^3)}{\mathcal{O}(r^3)\,\,\mathcal{O}(r^{1/2})}\rightarrow
0\,\, \textrm{as}\,\,r\rightarrow\infty\eqno(0.6) $$ This
original example of systolic freedom has been vastly
generalized by several authors (see \cite{BK} for an overview and
recent advances) to show that $\lq\lq$freedom" rather than
$\lq\lq$ridigidity" predominates for dimension $n\geq3$.

This left the case of $Z_2$ coefficients open for $n\geq3$. This
case has a remarkable relevance in quantum information theory,
which is the subject of another paper \cite{FLM}. Classically,
there is only one type of error: the $\lq\lq$bit flip."  In a
quantum mechanical context the algebra of possible errors has
two generators: $\lq\lq$bit flip" and $\lq\lq$relative phase."
It is possible to map the problem of correcting these (Fourier)
dual errors onto the problem of specifying (Poincare) dual
cycles in a manifold.  Torsion coefficients for the cycles
corresponds to finite dimensional quantum state spaces:
$Z_2$--coefficients correspond to expressing quantum states in
terms of qubits.

It is reported in \cite{KS} that Gromov conjectured $Z_2$--rigidity,
ie, systolic inequalities like $(0.1)$ and $(0.5)$ would hold
in this case of $Z_2$ coefficients.  The ease with which
nonoriented cycles can be modified to reduce area, particularly
in codimension equal to $1$, is well known in geometric measure
theory and lends support to the idea that at least
$Z_2-(n-1,1)$--rigidity might hold. In fact, the opposite is the
case. We will exhibit a family of Riemannian metrics on $S^2
\times S^1$ exhibiting
$Z_2-(2,1)$--systolic freedom: the ratios
$(Z_2-systole_{3}/Z_2-systole_{2} \cdot Z_{2}-systole_{1})$
approach zero as the parameter approaches infinity. Moreover,
from this example, as in \cite{BK}, quite general $Z_2$--freedoms
can be found.

In section 3, we discuss the quantification of systolic freedom
and note that the present example for $Z_2$--freedom is measured
by a function growing more slowly than $\log$ whereas in
Gromov's original example freedom grows by a power, and in an
example of Pittet \cite{P} freedom grows exponentially. It is now
of considerable interest, particularly in connection with
quantum information theory, whether the  $\lq\lq$weakness" of
$Z_2$--freedom is an artifact of the example or inherent.

\section{The Example}

As raw material, we use a sequence of closed hyperbolic surfaces
$\Sigma_g$ of genus $g\rightarrow\infty$ with the following
three properties:
\begin{enumerate}
\item[(i)]
$\lambda_1 (\Sigma_g) \geq c_1, \lambda_1$ being the smallest
eigenvalue of the Laplacian on functions,
\item[(ii)] There exits an isometry
$\tau\co \Sigma_g\rightarrow\Sigma_g$, with order $(\tau) \geq c_2
(\log\,g)^{1/2}$, and
\item[(iii)] The map $\Sigma_g \rightarrow
\Sigma_g/\tau(\sigma) \equiv \sigma$ is a covering projection
and the base surface
$_g{S}=:\Sigma_g/\tau(\sigma) \equiv \sigma$ has injectivity
radius $(_g{S}) \geq c_3(\log\,g)^{1/2}$
\end{enumerate} where $c_1,\, c_2,$ and $c_3$ are positive
constants independent of $g$.

We will return to the construction of the family $\{ \Sigma_g
\}$ at the end of this section. Let $M_g=(\Sigma_g\times
\mathbb{R})/(x,t) \equiv(\tau x,t+1)$ be the Riemannian
$\lq\lq$mapping torus" of $\tau$. We can also think of
$M_g=\Sigma_g \times [0,1]/(x,0)\equiv(\tau x,1)$. By two
theorems of Lickorish \cite{L}, we may first write $\tau^{-1}$ out
in the mapping class group of $\Sigma_g$ as a product of Dehn
twists
$\sigma_i$ along simple loops $\gamma_i \subset \Sigma_g$: $$
\tau^{-1} = \sigma_{n_g} \circ \dots \circ \sigma_2\, \circ\,
\sigma_1 \eqno{(1.1)} $$ and second perform Dehn surgeries along
pushed-in copies of $\{\gamma_i\}$: $$ \Big\{ \gamma_1
\times\Big(\frac{1}{2} + \frac{1}{3n_g}\Big),
\gamma_2\times\Big(\frac{1}{2}+\frac{2}{3n_g}\Big),\dots
\gamma_i
\times\Big(\frac{1}{2} + \frac{i}{3n_g}\Big),\dots \gamma_{n_g}
\times \Big(\frac{1}{2} +\frac{1}{3}\Big)\Big\} $$ to obtain a
diffeomorphic copy of $\Sigma_g \times[0,1]$ whose product
structure induces $[\tau^{-1}]\co \Sigma_g \times 0 \rightarrow
\Sigma_g \times 1$.

Thus, $n_g$ Dehn surgeries on $M_g$ produce the mapping torus
for
$\tau^{-1} \circ \tau$, ie $\Sigma_g \times S^1$.

In \cite{FLM}, we will find upper bounds on both $n_g$ and max
length
$(\gamma_i)$ in order to compute a lower bound on the
$Z_2$--freedom function. To merely establish $Z_2$--freedom, we
do not need these estimates. To convert $\Sigma_g \times S^1$
to $S^2 \times S^1$ an additional $2g$ Dehn surgeries are
needed: Do half (a $\lq\lq$sub kernel") of these surgeries at
level $\frac{1}{2} +\frac{1}{6n_g}$ and the dual half at level
$\frac{1}{2}$.  The result of all $n_g +2g$ Dehn surgeries is
topologically $S^2 \times S^1$, and once these surgeries are
metrically specified, we obtain a sequence of Riemannian
$3$--manifolds
$(S^2 \times S^1)_g=:S^{2} \times S_g^1$.

In section 2 where $Z_2$--freedom is established, four metrical
properties of these surgeries will be referenced.

They are:
\begin{itemize}
\item[(A)]
\emph{The core curves for the Dehn surgeries are taken, for
convenience, to be geodesics in $\Sigma_g \times[0,1]$ so that
the boundaries $\partial T_{i,\epsilon}$ of their $\epsilon$
neighborhoods are Euclidean flat.}\hfill(1.2)\break
\emph{Also $\epsilon > 0$ is chosen very small.  See (D).}
\item[(B)]
\emph{The replacement solid tori $T_{i,\epsilon}'$ have
$\partial T_{i,\epsilon}'$ isometric to $\partial
T_{i,\epsilon}$ and are defined as twisted products $D^2 \times
[0, 2 \Pi \epsilon]/
\beta$ where $\Pi(\epsilon)$ is a constant slightly larger than
$\pi$ so that the meridians in $T_{i,\epsilon}$ have length
$2\Pi\epsilon$ and $\beta$ is an isometric rotation of the disk
$D^2$ adjusted to equal the holonomy obtained by traveling
orthogonal to the surgery slopes in $\partial T_{i,\epsilon}$
from
$\partial D^{2} \times p t$ back to itself.}
\hfill(1.3)
\item[(C)]
\emph{The geometry on the disk $D^2$ above is rotationally
symmetric and has a product collar on its boundary as long as
the boundary itself.}\hfill(1.4)
\item[(D)]
\emph{Finally, $\epsilon > 0$ is so small that the total volume
of all the replacement solid tori,
$\cup_{i} T_{i,\epsilon}'$ is $o(g)$.}\hfill(1.5)
\end{itemize}

With specifications: $(A)\dots(D)$, Dehn surgery yields a
precise-smooth Riemannian manifold for which all the relevant
notions of $p$--area are defined.  We could work in this
category but there is no need to do so since perturbing to a
smooth metric will not effect the status of $(Z_2)$ systolic
freedom versus rigidity.

It is now time to return to the construction of the family
$\{\Sigma_g\}$.  We follow an approach of
\cite{SX} and \cite{Sc} in considering the co-compact torsion free
Fuchscian group $\Gamma_{(-1, p)}$, the group of unit norm
elements of the type $\frac{-1, p}{Q}$ quaternion algebra where
$p$ is prime and $p \equiv 3$ mod $4$. The group $\Gamma$ may
be explicitly written as:
$$
\Gamma_{(-1, p)}=
\bigg\{\left|a+b\sqrt{p} \quad -c+d\sqrt{p}
\atop c+d\sqrt{p} \quad a-b\sqrt{p}\right| : a, b, c, d \in Z,
det=1\bigg\} \biggl/ \pm \,\,\textrm{id.}\eqno{(1.6)}$$
Analogous to the congruence of $SL(2,R)$, we have for integers
$N>2$ the normal subgroups of $\Gamma_{-1,p}$,
$$
\Gamma_{(-1,p)}(N)=
\bigg\{\left|1+N(a+b \sqrt{p}) \quad N(-c+d\sqrt{p})
\atop  N(c+d\sqrt{p})\quad 1+N(a-b\sqrt{p})\right| : a, b, c, d
\in Z \bigg\}\biggl/ \pm\,\,\textrm{id.}
\eqno{(1.7)}$$ which are known (\cite{SX} and \cite{S}) to satisfy
(i).

In Lemma 2 \cite{Sc} it is proved that:
$$ inj. \,\, rad. \,(\mathbb{H}^2/\Gamma_{-1,p}(N))=
\mathcal{O}(\log N)
\eqno{(1.8)}$$ and in the proof of Theorem 6 that genus
$\bigl(\mathbb{H}^2/\Gamma_{-1,p} (N)\bigr) =: genus\,\,
\bigl(\Sigma(N)\bigr)=: genus\,(\Sigma_g)=:genus\,(N)$
satisfies:
$$
\mathcal{O}(N^2) \leq genus\,(N) \leq
\mathcal{O}(N^3)\eqno{(1.9)}
$$ so
$$ inj.\,\, rad.\, \Sigma_g = \mathcal{\mathcal{O}} (\log
g)\eqno{(1.10)}
$$

Now choose a sequence of $h$ and $g$ to satisfy $\log
g=\mathcal{O}(\log h)^2$ and so that $N(h)$ divides $N(g)$.
Thus, we have a covering projection
$\Sigma_g \longrightarrow \Sigma_h$.  Let $\alpha$ be the
shortest essential loop in $\Sigma_h$, by $(1.10)$
length$(\alpha)=\mathcal{O}(\log h)$. Choosing a base point on
$\alpha, \,[\alpha] \in
\Gamma_{(-1,p)} \bigl(N(h)\bigr)/\Gamma_{(-1, p)}
\bigl(N(g)\bigr)$ satisfies:
$$ order\, [\alpha] \geq
\mathcal{O}\bigl(\log(h)\bigr)=\mathcal{O} (\log g)^{1/2}
\eqno{(1.11)}
$$ since the translation length of $\alpha=\mathcal{O}(\log
g)^{1/2}$ must be multiplied by $\mathcal{O}(\log g)^{1/2}$
before it reaches length $\mathcal{O}(\log g)$, a necessary
condition to be an element in the subgroup
$\Gamma_{(-1,p)}\bigl(N(g)\bigr)$.

Let $\tau$ be the translation determined by $[\alpha]$.  We
have just checked condition (ii) order $(\tau)>\mathcal{O}(\log
g)^{1/2}$. Factor the previous covering as:
$$
\Sigma_g \longrightarrow \Sigma_g/\langle\tau\rangle
\longrightarrow
\Sigma_h \eqno{(1.12)}$$ and set $\Sigma_g/\langle\tau\rangle
=: {_g{S}}.$  Since $_g{S}$ covers
$\Sigma_h$, we conclude condition (iii):
$$ inj.\,\, rad.\,(_g{S})\geq inj.\,\, rad.\,(\Sigma_h)
\geq \mathcal{O} (\log h) = \mathcal{O} (\log
g)^{1/2}.\eqno{(1.13)}$$

\section{Verification of Freedom} 

We regard the Riemannian
manifold $S^2 \times S_g^1$ as essentially specified in
$\linebreak$ section 1.  Technically, there is the parameter
$\epsilon$ to be analyzed in \cite{FLM} which controls the
$\lq\lq$thickness" of the Dehn surgeries.  On two occasions, we
demand this to be sufficiently small (the cost is an increase
in the maximum absolute value of the Riemann curvature tensor
as a function of $g$). The first occurrence is in the next
proposition.
\begin{prop}
$vol\,(S^2 \times S_g^1) =Z_2-systole_3(S^2 \times
S_g^1)=\mathcal{O}(g)$
\end{prop}
\begin{proof}
$Volume\,(M_g)=vol\,\big( \Sigma_g \times[0,1]=area(\Sigma_g) =2
\pi \mathcal{X}(\Sigma_g)=\mathcal{O}(g) \big)$. By choosing
$\epsilon>0$ small enough as a function of $g$, the Dehn
fillings contribute negligible volume so this property is
retained by $S^2
\times S_g^1$.
\end{proof} The next proposition is more subtle.
\begin{prop}
$Z_2-systole_2(S^2 \times S_g^1)=
\mathcal{O}(g)$\end{prop}
\begin{proof} According to \cite{F} a  non-oriented minimizer
among all nonzero codimension one cycles always exists and is
smooth provided the ambient dimension is at most $7$.  Let $X_g
\subset S^2 \times S_g^1$ be this minimizer.  For a
contradiction, assume
$area(X_g)<\mathcal{O}(g)$.

The Dehn surgeries in section $1$ were confined to
$\Sigma_g \times [\frac{1}{2}, 1]$, so the surfaces
$\Sigma_g \times t$, $t\in(0,\frac{1}{2})$ persist as
submanifolds of $S^2 \times S_g^1$.  By Sard's theorem, for
almost all $t_\circ \in (0, \frac{1}{2})$,
$\Sigma_g \times t_\circ$ intersects $X_g$ transversely. Let
$W_t$, $t\in(0,\frac{1}{2})$ denote the intersection. By the
co-area formula.
$$
\mathcal{O}(g) > area (X_g) \geq
\int_{t=0}^{1/2}\,length\,(W_t)\,dt
\eqno{(2.1)}
$$

Consequently, for some transverse $t_\circ \in(0,\frac{1}{2})$,
$$ length\,(W_{t_\circ})< \mathcal{O}(g)
\eqno{(2.2)}
$$ Since both $\Sigma_g \times t_\circ$ and $X_g$ represent the
nonzero element of $H_2(S^2 \times S_g^1;Z_2)$, the complement
$S^2 \times S_g^1 \setminus (\Sigma_g \times t_\circ
\cup X_g)$ can be two colored into black and white regions
(change colors when crossing either surface) and the closure
$B$ of the black points is a piecewise smooth
$Z_2$--homology between $\Sigma_g \times t_\circ$ and $X_g$.

For homological reasons, the reverse Dehn surgeries
$S^2 \times S_g^1 \rightsquigarrow M_g$ have cores with zero
(mod 2) intersection with $X_g$.  This means that the tori
$\partial T_{i,\epsilon}=\partial T_{i,\epsilon}'$ each meet
$X_g$ in a null homologous, probably disconnected, $1$--manifold
$X_g \cap \partial T_{i,\epsilon}\, \subset \, \partial
T_{i,\epsilon}$. Again, if
$\epsilon$ is a sufficiently small function of $g$, we may
$\lq\lq$cut off"
$X_g$ along these tori to form
$$ X_{g}'=(X_g \setminus \cup_i T_{i,\epsilon})
\cup \delta_i,
$$ where $\delta_i$ denotes a bounding surface for
$X_g \cap \partial T_{i,\epsilon}$ in $\partial
T_{i,\epsilon}$, with negligible increase in area.  In
particular, we still have:
$$ area\,(X_{g}') < \mathcal{O}(g)
\eqno{(2.3)}
$$

More specifically choose $\delta_i$ to be the $\lq\lq$black"
piece of $\partial T _{i,\epsilon}$, ie $\delta_i \subset
B$.  If we set
$$B'=closure (B \setminus \cup_i T_{i,\epsilon})
$$ and recall
$$\cup_i T_{i,\epsilon} \cap \Sigma_g
\times t_\circ = \emptyset,
$$ we see that $B'$ is a $Z_2$--homology from $X_{g}'$ to
$\Sigma_g \times t_\circ$.

It is time to use property (i):  $W_{t_\circ}$ separates
$\Sigma_g \times t_\circ$ into two subsurfaces meeting along
their boundaries: One subsurface sees black on the positive
side, the other on its negative side.  An inequality of Buser's
\cite{Bu}, a converse to the Cheeger' isoperimetric inequality,
states that
$area>constant\cdot length$, in the presence of bounded
sectional curvatures, yields an upper bound on $\lambda_1$.
Thus, the smaller of these two subsurfaces, call it $Y \subset
\Sigma_g \times t_\circ$ must satisfy:
$$ area\,(Y)\leq c_4 \,\, length\, (W_{t_\circ})
\eqno{(2.4)}
$$ where $c_4$ is independent of $g$.  Combining with line
(2.2), we have:
$$ area \,(Y) \leq \mathcal{O}(\log g)
\eqno{(2.5)}
$$

Now modify $X_{g}'$ to $Z$ by cutting along $W_{t_\circ}$ and
inserting two parallel copies of $Y$.  This may be done so that
the result is disjoint from $\Sigma_g \times t_\circ$ but
bordant to it by a slight modification $B''$ of $B'$, with
$B''$ still disjoint from $\Sigma_g \times t_\circ$.  See
Figure 2.1 and Figure 2.2.

\begin{figure}[ht!]
\centerline{%
\relabelbox\small
\epsfxsize=120mm \epsfbox{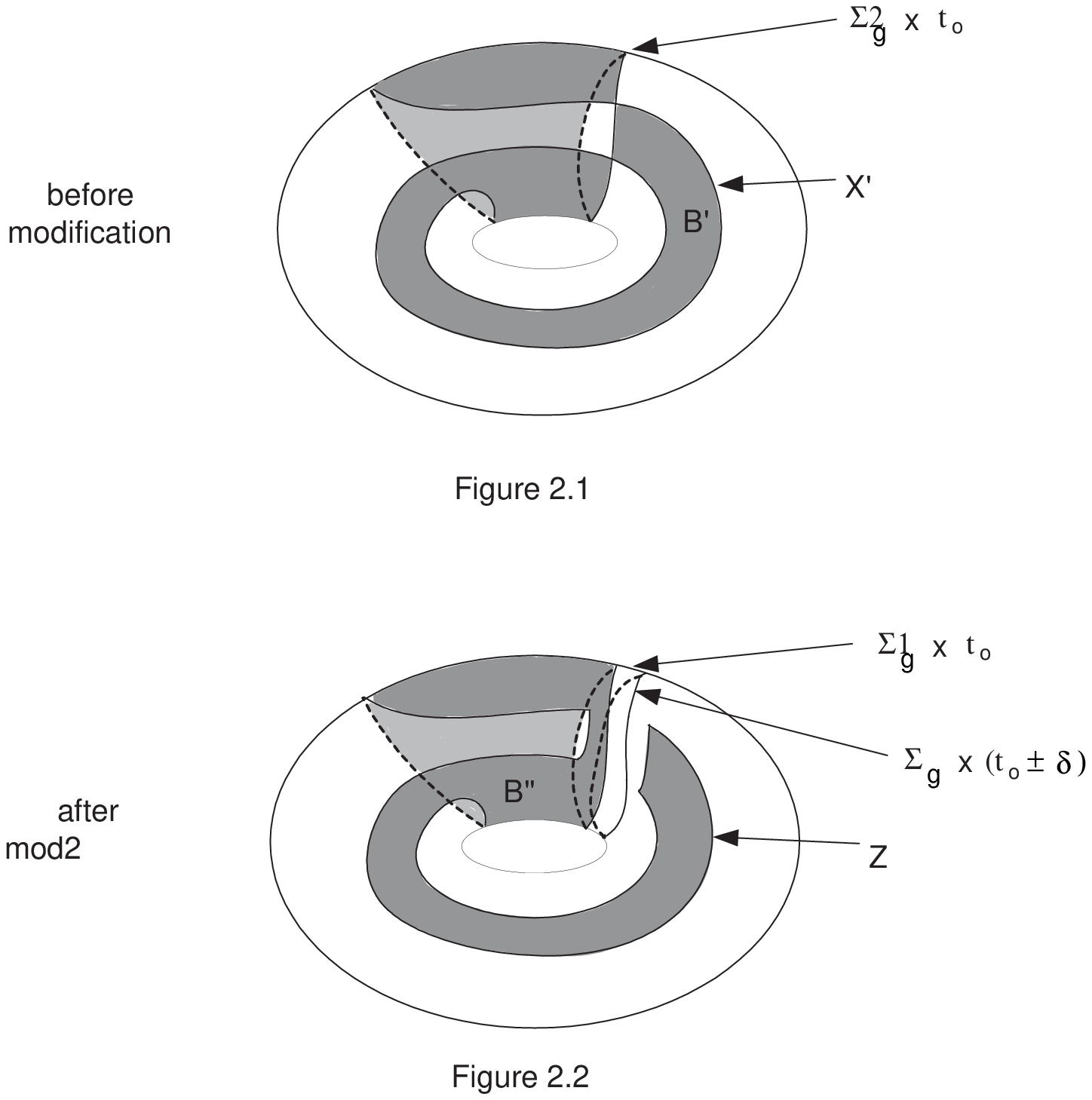}
\adjustrelabel <0pt,2pt> {Z}{$Z$}
\adjustrelabel <0pt,2pt> {X'}{$X'$}
\adjustrelabel <-4pt,0pt> {B'}{$B'$}
\adjustrelabel <0pt,2pt> {B''}{$B''$}
\relabel{S2}{$\Sigma_g\times t_0$}
\relabel{S1}{$\Sigma_g\times t_0$}
\relabel{S}{$\Sigma_g\times (t_0\pm \delta)$}
\relabel {Figure 2.1}{Figure 2.1}
\relabel {after}{after}
\relabel {before}{before}
\relabel {modification}{modification}
\relabel {mod2}{modification}
\relabel {Figure 2.2}{Figure\ \,2.2}
\endrelabelbox}
\end{figure}

combining (2.3) and (2.5):
$$ area\,(Z) \leq 3\cdot \mathcal{O}(\log g) =\mathcal{O}(\log
g)
\eqno{(2.6)}
$$

Now reverse the Dehn surgeries and consider:
$$ (B''; \, \Sigma_g \times t_\circ, Z)
\subset M_g \setminus \Sigma_g
\times (t_\circ \pm \delta) \subset M_{g}.
\eqno{(2.7)}
$$

The middle term of line (2.7) is diffeomorphic to
$\Sigma_g \times \mathbb{R}$, which is a codimension 0
submanifold of $\mathbb{R}^3$.  This proves that $B''$ and in
particular $Z$ is orientable.  But this looks absurd.
Apparently, we have constructed an oriented surface $Z$
oriented-homologous to the fiber $\Sigma_g \times t_\circ$ of
$M_g$ of smaller area (compare line (2.6) with the first line
in the proof of proposition 1.1).

Let $\frac{\partial}{\partial t}$ be the divergenceless flow in
the interval direction on $M_g$.  Lift $Z$ to $\widetilde{Z}$
in the infinite cyclic cover $\Sigma_g \times \mathbb{R}$ and
consider the flow through the lift $\widetilde{B}''$, the lift
of $B''$.  The divergence theorem states that the flux through
$\widetilde{Z}$ is equal to the flux through $\Sigma_g \times
t_\circ$.  Since $\frac{\partial} {\partial t}$ is orthogonal
to $\Sigma_g \times t_\circ$,
$$ area\,(\Sigma_g \times t_\circ) \leq area\,(\widetilde{Z}) =
area(Z)
\eqno{(2.8)}
$$ completing the contradiction.
\end{proof}
\begin{prop}
$Z_2-systole_{1}(S^2 \times S_g^1) \equiv
\mathcal{O} (\log g)^{1/2}$
\end{prop}
\begin{proof} We actually show that any homotopically essential
loop obeys this estimate.  The long collar condition $C$
(section 1) implies that any arc in $T'_{i,\epsilon}$ with end
points on $\partial T'_{i,\epsilon}$ can be replaced with a
shorter arc with the same end points lying entirely within
$\partial T'_{i,\epsilon}$.  It follows that any essential loop
in $S^2 \times S_g^1$ can be homotoped to a shorter loop lying
in the complement of the Dehn surgeries.

Thus, it is sufficient to show that any homotopically essential
loop
$\gamma$ in $M_g$ has $length\, \gamma \geq \mathcal{O} (\log
g)^{1/2}$.  For a contradiction, suppose the opposite.  Since
the bundle projection
$\pi \co  M_g \rightarrow [0,1]/\, 0\equiv 1$ is length
nonincreasing, degree
$\linebreak$
$\pi \circ \gamma < \mathcal{O} (\log g)^{1/2}$.  Lift $\gamma
\setminus pt$. to an arc
$\widetilde{\gamma}$ in $\Sigma_g \times R$.  The lift
$\widetilde{\gamma}$ joins some point $(p,t)$ to $(\tau^d p,
t+d)$ where $d=\,$ degree
$\pi \circ \gamma$. Since $d < \mathcal{O}(\log g)^{1/2}$ and
since condition (ii) requires order$(\tau)
\geq \mathcal{O}(\log g)^{1/2}$, we see that $p$ and $\tau^d p$
differ by a non-trivial  covering translation of the cover
$\Sigma_g \rightarrow {_g}S$.  Nevertheless, any non-trivial
covering translation moves each point of the total space at
least twice the injectivity radius of the base, a quantity
guaranteed by (iii) to be
$\geq \mathcal{O} (\log g)^{1/2}$.  Now using that the
projection
$\Sigma_g \times R \rightarrow \Sigma_g$ is also length
nonincreasing, we see that $length \,(\widetilde{\gamma}) \geq
\mathcal{O} (\log g)^{1/2}$. Since $length
\,(\widetilde{\gamma}) = length \,(\gamma)$, the same estimate
$\linebreak$ applies to $\gamma$.
\end{proof}

\begin{thm} The family $\{S^2 \times S_g^1\}$ exhibits
$Z_2$--systolic freedom.
\end{thm}

\proof From propositions 2.1, 2.2, and 2.3, we have:
$$\frac{Z_2-systole_3 (S^2 \times S_g^1)}{Z_2-systole_2 (S^2
\times S_g^1)
\cdot Z_2-systole_1 (S^2 \times S_g^1)} \leq
\frac{\mathcal{O}(g)}{\mathcal{O}(g)\,\,\,\mathcal{O}(\log
g)^{1/2} }\rightarrow 0.\eqno{\qed}$$

Many further examples in higher dimensions can now be generated.
It is easy to check that if $C$ is a circle of radius
$\frac{\mathcal{O}(g)}{\mathcal{O}(\log g)^{1/2}}$ then $(S^2
\times S_g^1)
\times C$ has
$\linebreak Z_2-(2,2)$--freedom.  As in \cite{BK}, two further
$1$--surgeries give a family of metrics on $S^2 \times S^2$ with
$Z_2-(2,2)$--freedom. Curiously, the homotopy theoretic methods
in \cite{BK} do not resolve whether
$CP^2$ has $Z_2$--freedom.  The difficulty is that a crucial
$\lq\lq$meromorphic map"  $CP^2 \rightarrow S^2 \times S^2$ has
even degree. Whether $CP^2$ admits a metric of
volume$\,=\epsilon$ in which every surface, orientable or not,
of area$\,\leq 1$ is null homotopic is an open question.  I
would like to thank M. Katz for his explanation of this
difficulty, and for orienting me within the literatures on
systolic inequalities.

\section{Curvature Normalization} The precise arithmetic of
both the theorem and Gromov's example (See introduction.)
suggests that the amount of systolic freedom exhibited in a
parameter family should be quantified.  The natural way to do
this is to homothetically rescale each metric in the family
(say $g$ is the parameter) to make the spaces as small as
possible while keeping all sectional curvatures bounded between
$-1$ and $+1$.

Given a family exhibiting $(p, q)$--freedom, for some choice of
coefficients, first rescale the members of the family to obtain
bounded curvature and then write the $\lq\lq
denominator"=systole_{p} (g) \cdot systole_{q} (g)$ as a
function of the rescaled $\lq\lq numerator"= volume=systole_{n}
(g)$. The function $F(n)=\frac{d(n)}{n}$ measures the
$\lq\lq$freeness" of the family.

In the constructions of Gromov and Babenko--Katz, $F(n)$ grows
like a positive power of $n$.  Pittet \cite{P} replaced a Nil
geometry construction of \cite{BK} with an analogous Solv geometry
construction to realize what our definition interprets as an
exponentially growing $F(n)$. When properly rescaled the growth
function for the examples in this paper will be considerably
slower than root $\log$ (to be estimated in \cite{FLM}). Perhaps
the most interesting question to arise from our example is
whether manifolds are $\lq\lq$nearly" $Z_2$--rigid, ie, do their
$Z_2$--freeness functions even when maximized over all families
of metrics grow with extreme slowness.  A negative answer would
be very interesting both within geometry and for the
implication for quantum codes. A positive answer would require
a new technical idea: eg, translating some as yet unproved
upper bound on the efficiency of quantum codes into
differential geometry.

I would like to thank David Meyer for may helpful discussions
on the theme of systolic geometry.

{\em The author is supported by Microsoft Research-Theory Group.}

\Addresses\recd
\end{document}